\documentclass[a4paper,12pt]{article}
\usepackage{amssymb}
\usepackage{amsmath}
\usepackage{amsthm}
\usepackage{graphicx}
\numberwithin{equation}{section}
\newcommand{\ep}{\varepsilon}
\newcommand{\la}{\lambda}
\newcommand{\va}{\varphi}
\newcommand{\ppp}{\partial}

\newcommand{\weight}{e^{2s\va}}

\newcommand{\R}{\mathbb{R}}

\newcommand{\p}{\mathbf{p}}
\newcommand{\www}{\widetilde}
\newcommand{\LLL}{\widetilde{L}_0}

\newcommand{\ssss}{\sum_{i,j=1}^n}
\newcommand{\oooo}{\overline}
\newcommand{\wwww}{\widetilde}

\allowdisplaybreaks

\title{Carleman estimate and unique continuation for 
a structured population model}

\author{Masaaki Uesaka $^1$ and Masahiro Yamamoto $^2$}

\date{}

\begin{document}
 \maketitle

\begin{center}
Graduate School of Mathematical Sciences, the University of Tokyo, 
3-8-1 Komaba, Meguro-ku, Tokyo, 153-8914, Japan.\\
E-mail: $^1$ muesaka@ms.u-tokyo.ac.jp $^2$ myama@ms.u-tokyo.ac.jp
\end{center}

\begin{abstract}
We consider a time-dependent 
structured population model equation and establish a
Carleman estimate.  We apply the Carleman estimate to prove the unique
continution which means that Cauchy data on any lateral 
boundary determine the solution uniquely.
\end{abstract}

\baselineskip 18pt
\section{Introduction}

Structured population models describe the change of distribution
of individuals in a population. In these models, individuals are described
by using several parameters-- for example, age, size and so on-- and
a population density is considered as a function of not only time
and spatial position but these individual parameters. In this meaning,
we can say that structured population models describe {}``the detail
structure of population.'' These models originated in Sharpe and
Lotka \cite{ShaLot} and McKendrik \cite{McK} and have been widely
studied in the mathematical biology. 

In this paper, we consider one of structured population models stated
in Webb \cite{Web} in which age and size are considered as individual
parameters. The model is described as follows: Let $\Omega\subset
\mathbb{R}^n$
be an open set which represents an inhabited area and $a_{1},\tau_{1},\tau_{2},
T$ be positive real constants. 
Henceforth $g>0$ on $[\tau_1,\tau_2]$ and $g\in C^1[\tau_1,\tau_2]$, and
we set
$$
Kv(x,t,a,\tau) = \sum_{i,j=1}^n a_{ij}(x)\ppp_i\ppp_jv(x,t,a,\tau) 
- \sum_{k=1}^n b_k(x,t,a,\tau)\ppp_kv - c(x,t,a,\tau)v, \quad
x \in \Omega, \thinspace 0 < t < T
$$
where $a_{ij} = a_{ji} \in C^1(\overline{\Omega})$, $1 \le i,j\le n$,
$b_k, c\in L^{\infty}(\Omega \times (0,T)\times (0,a_1)\times 
(\tau_1,\tau_2))$, and there exists 
a constant $\sigma_1>0$ such that 
$$
\sum_{i,j=1}^n a_{ij}(x)\xi_i\xi_j \ge \sigma_1 \sum_{i=1}^n\xi_i^2
$$
for $x \in \overline{\Omega}$ and $\xi_1, ..., \xi_n \in \R$.

Then our model equation is
$$
\ppp_t u(x,t,a,\tau) + \ppp_au(x,t,a,\tau) + \ppp_\tau(g(\tau)u(x,t,a,\tau))
= Ku(x,t,a,\tau),
$$
$$
\qquad \qquad 
(x,t,a,\tau)\in \Omega\times (0,T)\times(0,a_{1})\times(\tau_{1},\tau_{2}),
                                             \eqno{(1.1)}
$$
with initial and boundary conditions
$$
u(x,t,0,\tau)=\int_{0}^{a_{1}}\int_{\tau_{1}}^{\tau_{2}}
\beta\left(x,a,\tau,\widetilde{\tau}\right)u(x,t,a,\widetilde{\tau}) 
d\widetilde{\tau}da,
$$
$$
\qquad \qquad \qquad \qquad (x,t,\tau)\in(0,T)\times\Omega\times
(\tau_{1}, \tau_{2}),
                                           \eqno{(1.2)}
$$
$$
u(x,t,a,\tau_{1})=0, \quad(x,t,a)\in \Omega \times (0,T) \times
(0,a_{1}),                                    \eqno{(1.3)}
$$
$$
u(x,0,a,\tau)=p(x,a,\tau), \quad(x,a,\tau)\in\Omega\times(0,a_{1})
\times(\tau_{1},\tau_{2}),
                                            \eqno{(1.4)}
$$
and
$$
\partial_{\nu}u=0\quad\mbox{on $\ppp\Omega \times (0,T)
\times(0,a_{1})\times(\tau_{1},\tau_{2})$}.               \eqno{(1.5)}
$$

We can interpret equation (1.1) as follows.  The variable $a$ is the 
age of individual and $\tau$ the size, $u(x,t,a,\tau)$ can be interpreted
as the population density at time $t$, position $x$, age $a$
and size $\tau$. Moreover
\begin{itemize}
\item $\ppp_a u(x,t,a,\tau)$ represents \emph{aging effect}. 
The coefficient
is always exactly 1 because age increases exactly 1 per a year.
\item 
$\partial_{\tau}\left(g(\tau)u(t,x,a,\tau)\right)$ represents \emph{growth
effect }with $g(\tau)$ a \emph{growth modulus }, that is, 
$\int_{\tau}^{\tau'}1/g(\sigma) d\sigma$
is a spending time to grow the individual from size $\tau$ to size 
$\tau'$).
\item 
The ellipic part $Ku$ represents diffusion and taxis.  In particular,
$c(x,t,a,\tau)$ denotes a linearized mortality rate.
\item 
The condition (1.2) represents birth with birth rate $\beta$.
\end{itemize}
For details, see Webb \cite{Web} which also proves 
the existence of the solution of the system
(1.1) - (1.5) by the semigroup theory.

In this paper, we discuss 
\\
{\bf Unique continuation}.  Let $\Gamma \subset \ppp\Omega$ be
an arbitrary subboundary.  Then determine $u$ in 
$\Omega \times (0,T)\times (0,a_1) \times (\tau_1,\tau_2)$ by 
lateral Cauchy data $u, \nabla u$ on 
$\Gamma \times (0,T)\times (0,a_1) \times (\tau_1,\tau_2)$.
\\

For the unique continuation for elliptic, parabolic and 
hyperbolic equations, there are many works and see for example
Bellassoued and Yamamoto \cite{BY}, H\"ormander \cite{Ho}, 
Isakov \cite{I}, Lavrent'ev, Romanov and Shishat$\cdot$ski\u\i
\cite{LRS} and the references therein. 
Here we do not intend any comprehensive lists of references.
However as for the unique continuation for equation (1.1), 
to the authors' best knowledge, there are no results published.

We state our main result:
\\
{\bf Theorem 1 (unique continuation)}. 
Let $\Gamma \subset\ppp\Omega$ be an arbitrary subboundary.
Let $u \in H^1(\Omega \times (0,T)\times (0,a_1) \times (\tau_1,\tau_2))$ 
satisfy $\ppp_i\ppp_ju \in L^2(\Omega \times (0,T)\times (0,a_1) \times
(\tau_1,\tau_2))$ and (1.1).
Then $u = \vert \nabla u\vert = 0$ on 
$\Gamma \times (0,T)\times (0,a_1) \times (\tau_1,\tau_2)$ yields  
$u=0$ in $\Omega \times (0,T)\times (0,a_1) \times (\tau_1,\tau_2)$.
\\

The proof is based on a Carleman estimate for (1.1), which may be 
interesting itself.
In the case of $g\equiv 0$, we refer to Traore \cite{Tr}, which proves
a Carleman estimate with a weight function in the form of
$\exp\left(2s\frac{V(x)}{at(T-t)}\right)$ with some function $V(x)$.
The weight function in \cite{Tr} is inspired by \cite{ImaF} and \cite{Ima},
and \cite{Tr} discusses the controllability, but it it very difficult
to derive the unique continuation by the Carleman estimate 
in \cite{Tr}.
\\

The paper is composed of three sections.  In Section 2, we prove the 
key Carleman estimate for (1.1) and in Section 3, the proof of 
Theorem 1 is completed.
\section{Carleman estimate}

We set 
$$
L_0u := \ppp_tu + \ppp_au + \ppp_{\tau}(g(\tau)u).
$$
Let $D \subset \Omega \times (0,T)\times (0,a_1) \times (\tau_1,\tau_2)$ be 
a subdomain.
\\
{\bf Lemma 1}. Let $d \in C^2(\overline{\Omega})$ satisfy 
$\vert \nabla d(x)\vert \ne 0$ on $\overline{\Omega}$.
We fix $t_0 \in (0,T)$, $a_0 \in (0,a_1)$ and $\tau_0 \in (\tau_1,\tau_2)$ 
arbitrarily and set
$$
\psi(x,t,a,\tau) = d(x) - \beta(\vert t-t_0\vert^2
+ \vert a-a_0\vert^2+\vert \tau-\tau_0\vert^2)
$$
and
$$
\va(x,t,a,\tau) = e^{\lambda\psi(x,t,a,\tau)}.
$$
Then there exists a constant $\la_0>0$ 
such that for arbitrary $\la \ge \la_0$,
we can choose a constant $s_0(\la) > 0$ satisfying: there exists 
a constant $C=C(s_0,\la_0)>0$ such that
$$
\int_D\left\{ \frac{1}{s\va}\vert L_0u\vert^2
+ s\la^2\va\vert \nabla u\vert^{2} + s^{3}\la^4\va^3 u^{2}\right\}
\weight dxdtdad\tau\\
\leq C\int_D \vert (L_0-K)u\vert^2\weight dxdtdad\tau             \eqno{(2.1)}
$$
for all $s > s_0$ and all $u\in H^1(D)$ satisfying
$\ppp_i\ppp_ju \in L^2(D)$ and $\text{supp} \thinspace u \in D$.
\\

The constant $C > 0$ in (2.1) depends continuously on
$$
\max_{1\le i,j\le n}\Vert a_{ij}\Vert_{C^1(\overline{D})},\quad
\Vert b_i\Vert_{L^{\infty}(D)}, \quad 
\Vert c \Vert_{L^{\infty}(D)}.
$$

{\bf Remark.}  We further assume that for each
$\www{t}, \www{a}, \www{\tau}$, the cross section 
$\{ x; \thinspace (x,\www{t}, \www{a}, \www{\tau}) \in D\}$
is composed of a finite number of smooth surfaces.  
Then similarly to Theorem 3.1 in Yamamoto \cite{Y}, we can 
improve (2.1) as 
\begin{align*}
& \int_D\Biggl\{ \frac{1}{s\va}\left(
\vert L_0u\vert^2
+ \sum_{i,j=1}^n \vert \ppp_i\ppp_ju\vert^2\right)
+ s\la^2\va\vert \nabla u\vert^{2} + s^{3}\la^4\va^3 u^{2}\Biggr\}
\weight dxdtdad\tau\\
\leq &C\int_D \vert (L_0-K)u\vert^2\weight dxdtdad\tau.
\end{align*}
\vspace{0.3cm}
\\
{\bf Proof of Lemma 1}.   
We set 
$$
\www{L}_0u := \ppp_tu + \ppp_au + g(\tau)\ppp_{\tau}u.
$$
Then $L_0u = \www{L}_0u + g'(\tau)u$.
Thanks to the large parameters $s$ and
$\lambda$, it is sufficient to prove the Carleman estimate for
$$
Lu:= \LLL u - \sum_{i,j=1}^n a_{ij}\ppp_i\ppp_ju.
$$
By a usual density argument (i.e., the approximation of any $u$
satisfying the condition in Lemma 1 by a sequence 
$u_m \in C^{\infty}_0(D)$), it suffices to prove  
the Carleman estimate for $u\in C^{\infty}_0(D)$.

We further set
$$
\sigma(x) = \ssss a_{ij}(x)(\ppp_id)(x)(\ppp_jd)(x), \quad
x \in \overline{\Omega}
$$
and
$$
w(x,t) = e^{s\va(x,t)}u(x,t)
$$
and
$$
Pw(x,t) = e^{s\va}L(e^{-s\va}w) = e^{s\va}Lu.   \eqno{(2.2)}
$$

The proof is very similar to the proof of Theorem 3.1 in 
Yamamoto \cite{Y}, which is a Carleman 
estimate for a parabolic equation.  More precisely,
\begin{enumerate}
\item
the decomposition of $P$ into the part $P_1$ and $P_2$, where 
$P_1$ is composed of second-order and zeroth-order 
terms in $x$, and $P_2$ is composed of first-order terms in $t$ 
and first-order terms in $x$.  
\item
Estimation of $\int_D (\vert P_2w\vert^2 + 2(P_1 w)(P_2w))dxdtdad\tau$
from the below.
\item
Another estimate for 
$$
\int_D Pw \times \text{[the term $u$ with
second highest order of $s, \la, \va$ among $Pw$]}.
$$
\end{enumerate}

Direct calculation of (2.2) gives 
$$
Pw = \LLL w - \sum_{i,j=1}^n a_{ij}(x)\ppp_i\ppp_jw
+ 2s\la\va \ssss a_{ij}(x)(\ppp_id)\ppp_jw
$$
$$
- s^2\la^2\va^2 \sigma w + s\la^2\va \sigma w 
+ s\la\va w \ssss a_{ij}\ppp_i\ppp_jd - s\la\va w(\LLL\psi)\thinspace
\mbox{in $D$}.                          \eqno{(2.3)}
$$
Here we note that we have specified all the dependency of 
coefficients on $s$, $\la$ and $\va$.
We set
$$
A_1 = s\la^2\va\sigma + s\la\va \ssss a_{ij}\ppp_i\ppp_jd
- s\la\va(\LLL\psi)
=: s\la^2\va a_1(x,t,a,\tau;s,\la).
$$
Then
\begin{align*}
&Pw = \LLL w - \sum_{i,j=1}^n a_{ij}(x)\ppp_i\ppp_jw
+ 2s\la\va \ssss a_{ij}(x)(\ppp_id)\ppp_jw \\
-& s^2\la^2\va^2 \sigma w + A_1w = fe^{s\va}\quad \text{in $D$}.
\end{align*}
We note that $a_1$ depends on $s$ and $\la$ but  
$$
\vert a_1(x,t,a,\tau;s,\la)\vert \le C
$$
for $(x,t,a,\tau) \in \oooo D$ and all sufficiently large $\la > 0$ and 
$s > 0$.
Here and henceforth by $C$, $C_1$, etc., we denote generic constants
which are independent of $s$, $\la$ and $\va$ but may change 
line by line.
 
Then taking into consideration the orders of $(s,\la,\va)$, 
we divide $Pw$ as follows:
$$
P_1w = - \sum_{i,j=1}^n a_{ij}(x)\ppp_i\ppp_jw
- s^2\la^2\va^2 w \sigma(x,t) + A_1w                 \eqno{(2.4)}
$$
and
$$
P_2w = \LLL w + 2s\la\va \ssss a_{ij}(x)(\ppp_id)\ppp_jw.  \eqno{(2.5)}
$$
By $\Vert fe^{s\va}\Vert^2_{L^2(D)}
= \Vert P_1w + P_2w\Vert^2_{L^2(D)}$, we have
$$
2\int_D (P_1w)(P_2w) dxdtdad\tau + \Vert P_2w\Vert_{L^2(D)}^2
\le \int_D f^2\weight dxdtdad\tau.                               \eqno{(2.6)}
$$
We estimate:
\begin{align*}
&\int_D (P_1w)(P_2w) dxdtdad\tau 
= - \ssss \int_D a_{ij}(\ppp_i\ppp_jw) (\LLL w) dxdtdad\tau\\
- &\ssss \int_D a_{ij}(\ppp_i\ppp_jw) 2s\la\va
\sum_{k,\ell=1}^n a_{k\ell}(\ppp_kd)(\ppp_{\ell}w) dxdtdad\tau\\
- & \int_D s^2\la^2\va^2\sigma w(\LLL w) dxdtdad\tau
- \int_D 2s^3\la^3\va^3\sigma w \ssss a_{ij}(\ppp_id)(\ppp_jw) dxdtdad\tau\\
+ & \int_D (A_1w)(\LLL w) dxdtdad\tau + \int_D (A_1w)2s\la\va \ssss
a_{ij}(\ppp_id)(\ppp_jw) dxdtdad\tau
\end{align*}
$$
=: \sum_{k=1}^6 J_k.                       \eqno{(2.7)}
$$
Now, applying the integration by parts, $a_{ij} = a_{ji}$
and $u \in C^{\infty}_0(D)$ and assuming that $\la > 1$ and $s > 1$ 
are sufficiently large, we reduce all the derivatives of $w$ to
$w, \ppp_iw, \LLL w$.  
We note that 
$$
\int_D u(\LLL v) dxdtdad\tau = -\int_D (L_0u)v dxdtdad\tau,
\quad u,v \in C^1_0(D).
$$
We continue the estimation of $J_k$, $1 \le k \le 6$.
\begin{align*}
&\vert J_1\vert = \left\vert -\ssss \int_D 
a_{ij}(\ppp_i\ppp_jw)(\LLL w) dxdtdad\tau \right\vert \\
= &\left\vert \ssss \int_D (\ppp_ia_{ij})(\ppp_jw)(\LLL w) dxdtdad\tau 
+ \ssss \int_D a_{ij} (\ppp_jw)\ppp_i(\LLL w) dxdtdad\tau\right\vert\\
=& \Biggl\vert \ssss \int_D (\ppp_ia_{ij})(\ppp_jw)(\LLL w) dxdtdad\tau \\
+& \Biggl(\sum_{i>j} \int_D a_{ij}((\ppp_jw)\ppp_i(\LLL w) 
+ (\ppp_iw)\ppp_j(\LLL w)) dxdtdad\tau \\
+ &\int_D \sum_{i=1}^n a_{ii}(\ppp_iw)\ppp_i(\LLL w) dxdtdad\tau\Biggr)
\Biggr\vert 
\end{align*}
$$
\le C\int_D \vert \nabla w\vert \vert \LLL w\vert dxdtdad\tau.
                                                \eqno{(2.8)}
$$

Here we used
$$
\LLL ((\ppp_iw)\ppp_jw) = (\LLL(\ppp_iw))\ppp_jw
+ (\ppp_iw)\LLL(\ppp_jw),
$$
and
\begin{align*}
& \Biggl(\sum_{i>j} \int_D a_{ij}((\ppp_jw)(\ppp_i\LLL w) 
+ (\ppp_iw)(\ppp_j\LLL w)) dxdtdad\tau \\
+ &\int_D \sum_{i=1}^n a_{ii}(\ppp_iw)(\ppp_i\LLL w) dxdtdad\tau\Biggr)
= \frac{1}{2} \ssss \int_D a_{ij}\LLL((\ppp_jw)(\ppp_iw)) dxdtdad\tau\\
= &-\frac{1}{2}\ssss \int_D L_0(a_{ij})(\ppp_jw)(\ppp_iw) dxdtdad\tau
= 0,
\end{align*} 
because $a_{ij}$ are independent of $t,a,\tau$.

Next 
\begin{align*}
&J_2 = -\ssss \sum_{k,\ell=1}^n \int_D 2s\la\va a_{ij}a_{k\ell}(\ppp_kd)
(\ppp_{\ell}w)(\ppp_i\ppp_jw) dxdtdad\tau\\
=& 2s\la \int_D  \ssss \sum_{k,\ell=1}^n \la(\ppp_id)\va
a_{ij}a_{k\ell}(\ppp_kd)(\ppp_{\ell}w)(\ppp_jw)dxdtdad\tau\\
+ &2s\la\int_D \ssss \sum_{k,\ell=1}^n \va \ppp_i(a_{ij}a_{k\ell}
\ppp_kd)(\ppp_{\ell}w)(\ppp_iw) dxdtdad\tau\\
+& 2s\la \int_D \ssss \sum_{k,\ell=1}^n \va a_{ij}a_{k\ell}
(\ppp_kd)(\ppp_i\ppp_{\ell}w)(\ppp_jw) dxdtdad\tau.
\end{align*}
We have 
$$
\mbox{[first term]}
= 2s\la^2\int_D \va \left\vert \ssss a_{ij}(\ppp_id)(\ppp_jw)
\right\vert^2 dxdtdad\tau \ge 0,
$$
and similarly to $J_1$, we can estimate
\begin{align*}
&\text{[third term]}
= s\la\ssss \sum_{k,\ell=1}^n \int_D 
\va a_{ij}a_{k\ell}(\ppp_kd)\ppp_{\ell}((\ppp_iw)(\ppp_jw))\\
= & - s\la^2\int_D \va\sigma \ssss a_{ij}(\ppp_iw)(\ppp_jw) dxdtdad\tau
- s\la \int_D \va \ssss \sum_{k,\ell=1}^n  
\ppp_{\ell}(a_{ij}a_{k\ell}\ppp_kd) (\ppp_iw)(\ppp_jw) dxdtdad\tau.
\end{align*}
Hence
\begin{align*}
&J_2
\ge - \int_D s\la^2\va\sigma \ssss a_{ij}(\ppp_iw)(\ppp_jw) dxdtdad\tau\\
- & C\int_D s\la\va \vert \nabla w\vert^2 dxdtdad\tau
+ 2s\la^2\int_D \va \left\vert \ssss a_{ij}(\ppp_id)(\ppp_jw)
\right\vert^2 dxdtdad\tau
\end{align*}
$$
\ge - \int_D s\la^2\va\sigma\ssss a_{ij}(\ppp_iw)(\ppp_jw) dxdtdad\tau
- C\int_D s\la\va \vert \nabla w\vert^2 dxdtdad\tau.               \eqno{(2.9)}
$$
$$
\vert J_3 \vert = \left\vert
- \int_D \frac{1}{2} s^2\la^2\va^2\sigma \LLL(w^2)dxdtdad\tau \right\vert
= \frac{1}{2}\left\vert \int_D L_0(s^2\la^2\va^2\sigma)w^2 dxdtdad\tau
\right\vert
$$
$$
\le C\int_D s^2\la^3\va^2 w^2 dxdtdad\tau.                \eqno{(2.10)}
$$
\begin{align*}
&J_4 = - \int_D 2s^3\la^3\va^3\sigma w \ssss a_{ij}(\ppp_id)
(\ppp_jw) dxdtdad\tau\\
= &- \int_D s^3\la^3\va^3 \ssss \sigma a_{ij}(\ppp_id)\ppp_j(w^2) dxdtdad\tau\\
=& \int_D s^3\la^3 \ssss 3\va^2\{ \la(\ppp_jd)\va\}
\sigma a_{ij}(\ppp_id)w^2 dxdtdad\tau\\
+ &\int_D s^3\la^3\va^3 \ssss \ppp_j(\sigma a_{ij}\ppp_id)w^2
dxdtdad\tau 
\end{align*}
$$
\ge \int_D 3s^3\la^4\va^3\sigma^2 w^2 dxdtdad\tau 
- C\int_D s^3\la^3\va^3 w^2 dxdtdad\tau.                          \eqno{(2.11)}
$$
\begin{align*}
&\vert J_5\vert
= \left\vert \int_D (A_1w)(\LLL w) dxdtdad\tau \right\vert
= \left\vert \int_D s\la^2\va a_1 w(\LLL w) dxdtdad\tau\right\vert\\
= &\frac{1}{2}\left\vert \int_D s\la^2\va a_1\LLL(w^2) dxdtdad\tau 
\right\vert
= \frac{1}{2}\left\vert \int_D s\la^2L_0(\va a_1)w^2 dxdtdad\tau 
\right\vert
\end{align*}
$$
\le  C\int_D s\la^3\va w^2 dxdtdad\tau.                  \eqno{(2.12)}
$$
\begin{align*}
&\vert J_6\vert  
= \left\vert \int_D s\la^2\va a_1 \times 2s\la\va w \ssss a_{ij}
(\ppp_id)(\ppp_jw) dxdtdad\tau \right\vert \\
= &\left\vert \int_D 2a_1s^2\la^3\va^2 
\ssss a_{ij}(\ppp_id)w(\ppp_jw) dxdtdad\tau \right\vert\\
= & \left\vert \int_D a_1s^2\la^3\va^2 
\ssss a_{ij}(\ppp_id)\ppp_j(w^2) dxdtdad\tau \right\vert\\
=& \left\vert -\int_D \sum_{i,j=1}^n\ppp_j(a_1s^2\la^3\va^2a_{ij}(\ppp_id))
w^2 dxdtdad\tau \right\vert
\end{align*}
$$ 
\le  C \int_D s^2\la^4\va^2 w^2 dxdtdad\tau.                 \eqno{(2.13)}
$$

Hence, choosing $s>0$ and $\lambda>0$ large, by (2.7) - (2.13) we obtain
\begin{align*}
& \int_D (P_1w)(P_2w) dxdtdad\tau 
\ge 3\int_D s^3\la^4\va^3\sigma^2 w^2 dxdtdad\tau
- \int_D s\la^2\va\sigma \ssss a_{ij}(\ppp_iw)(\ppp_jw) dxdtdad\tau\\
-& C\int_D s\la\va \vert \nabla w\vert^2 dxdtdad\tau
- C\int_D (s^3\la^3\va^3 + s^2\la^4\va^2)w^2 dxdtdad\tau
- C\int_D \vert \nabla w\vert \vert \LLL w\vert dxdtdad\tau.   
\end{align*}

Consequently
\begin{align*}
&3\int_D s^3\la^4\va^3\sigma^2 w^2 dxdtdad\tau
- \int_D s\la^2\va\sigma \ssss a_{ij}(\ppp_iw)(\ppp_jw) dxdtdad\tau\\
\le & \int_D (P_1w)(P_2w) dxdtdad\tau 
+ C\int_D s\la\va \vert \nabla w\vert^2 dxdtdad\tau
\end{align*}
$$
+ C\int_D (s^3\la^3\va^3 + s^2\la^4\va^2)w^2 dxdtdad\tau
+ C\int_D \vert \nabla w\vert \vert \LLL w\vert dxdtdad\tau.
                                                \eqno{(2.14)}
$$

Moreover for all large $s > 0$, by the definition (2.5) of $P_2$ and
an inequality: \\
$\vert \alpha + \beta\vert^2 \ge \frac{1}{2}
\vert \alpha\vert^2 - \vert \beta\vert^2$, we obtain
\begin{align*}
& \int_D \vert P_2w\vert^2 dxdtdad\tau
\ge \int_D \frac{1}{s\va} \vert P_2w\vert^2 dxdtdad\tau\\
= &\int_D \frac{1}{s\va}\left\vert \LLL w
+ 2s\la\va \ssss a_{ij}(\ppp_id)(\ppp_jw)\right\vert^2 dxdtdad\tau
\end{align*}
$$
\ge \frac{1}{2}\int_D \frac{1}{s\va}\vert \LLL w\vert^2 dxdtdad\tau
- C\int_D s\la^2\va \left\vert \ssss a_{ij}
(\ppp_id)(\ppp_jw)\right\vert^2 dxdtdad\tau,
$$
that is,
$$
\ep\int_D \frac{1}{s\va}\vert \LLL w\vert^2 dxdtdad\tau
\le C\ep\int_D \vert P_2w\vert^2 dxdtdad\tau
+ C\ep\int_D s\la^2\va \vert \nabla w\vert^2 dxdtdad\tau
$$
for any $\ep > 0$.
Hence by (2.14) and (2.6), we have
\begin{align*}
&3\int_D s^3\la^4\va^3\sigma^2 w^2 dxdtdad\tau
- \int_D s\la^2\va\sigma \ssss a_{ij}(\ppp_iw)(\ppp_jw) dxdtdad\tau\\
+ & \ep\int_D \frac{1}{s\va} \vert \LLL w \vert^2 dxdtdad\tau \\
\le & C\int_D f^2\weight dxdtdad\tau 
+ C\int_D s\la\va \vert \nabla w\vert^2 dxdtdad\tau
+ C\ep\int_D s\la^2\va \vert \nabla w\vert^2 dxdtdad\tau \\
+ &C\int_D (s^3\la^3\va^3 + s^2\la^4\va^2)w^2 dxdtdad\tau
+ C\int_D \vert \nabla w\vert \vert \LLL w\vert dxdtdad\tau.
\end{align*}

Now we note that the factor with the maximal order in $s, \la, \va$
of $w^2$ is $s^3\la^4\va^3\sigma^2$, the maximal factor of 
$\vert \nabla w\vert^2$ is $s\la^2\va\sigma$, and the maximal order 
of $\vert \LLL w\vert^2$ is $\frac{1}{s\va}$.  For example,
since we can choose $s, \la$ large, the term 
$(s^3\la^3\va^3 + s^2\la^4\va^2)w^2$ is of lower order.

Here, since the Cauchy-Schwarz inequality implies
\begin{align*}
&\vert \LLL w\vert \vert \nabla w\vert 
= s^{-\frac{1}{2}}\va^{-\frac{1}{2}}\la^{-\frac{1}{2}}
\vert \LLL w \vert
s^{\frac{1}{2}}\va^{\frac{1}{2}}\la^{\frac{1}{2}}\vert\nabla w \vert\\
\le & \frac{1}{2}\frac{1}{s\la\va}\vert \LLL w \vert^2
+ \frac{1}{2}s\la\va \vert \nabla w \vert^2,
\end{align*}
we have
\begin{align*}
&3\int_D s^3\la^4\va^3\sigma^2 w^2 dxdtdad\tau
- \int_D s\la^2\va\sigma \ssss a_{ij}(\ppp_iw)(\ppp_jw) dxdtdad\tau\\
+ & \left( \ep-\frac{C}{\la}\right)\int_D \frac{1}{s\va} 
\vert \LLL w \vert^2 dxdtdad\tau \\
\le & C\int_D f^2\weight dxdtdad\tau 
+ C\int_D s\la\va \vert \nabla w\vert^2 dxdtdad\tau
+ C\ep\int_D s\la^2\va \vert \nabla w\vert^2 dxdtdad\tau 
\end{align*}
$$
+ C\int_D (s^3\la^3\va^3 + s^2\la^4\va^2)w^2 dxdtdad\tau.    \eqno{(2.15)}
$$

The first and the second terms on the left-hand side of (2.15) have 
different signs and so we need another estimate.
Thus we will execute another estimation for
$$
\int_D s\la^2\va \sigma \ssss a_{ij}(\ppp_iw)(\ppp_jw) dxdtdad\tau
$$
by means of 
$$
\int_D (P_1w+P_2w) \times (s\la^2\va\sigma w) dxdtdad\tau.     
$$
Here we have chosen the factor $s\la^2\va\sigma w$ for obtaining
the term of $\vert \nabla w\vert^2$ with desirable $(s,\la,\va)$-factor
$s\la^2\va$.  That is, multiplying
$$
\LLL w + 2s\la\va \ssss a_{ij}(\ppp_id)(\ppp_jw)
- \ssss a_{ij}\ppp_i\ppp_jw
- s^2\la^2\va^2\sigma w + A_1w = fe^{s\va}
$$
with $s\la^2\va \sigma w$, we have
\begin{align*}
&\int_D (\LLL w)(s\la^2\va \sigma w)dxdtdad\tau 
+ \int_D 2s\la\va\ssss a_{ij}(\ppp_id)(\ppp_jw)s\la^2\va \sigma w
dxdtdad\tau \\
-& \int_D \left( \ssss a_{ij}\ppp_i\ppp_jw\right) s\la^2\va \sigma w
dxdtdad\tau - \int_D s^3\la^4\va^3\sigma^2 w^2 dxdtdad\tau\\
+ &\int_D (A_1w)(s\la^2\va \sigma w) dxdtdad\tau
\end{align*}
$$
=: \sum_{k=1}^5 I_k = \int_D fe^{s\va}s\la^2\va \sigma w dxdtdad\tau.
                                                 \eqno{(2.16)}
$$
Now, in terms of the integration by parts and $w \in C^2_0(D)$,
noting that $\vert \LLL\va\vert = \vert \la(\LLL\psi)\va\vert
\le C\la\va$ and $\ppp_i\va = \la(\ppp_id)\va$, etc.,
we estimate the terms.
$$
\vert I_1\vert = \left\vert \int_D \frac{1}{2}s\la^2\va\sigma 
\LLL(w^2) dxdtdad\tau \right\vert
\le C\int_D s\la^3\va w^2 dxdtdad\tau.              \eqno{(2.17)}
$$
\begin{align*}
&\vert I_2 \vert 
= \left\vert \int_D s^2\la^3\va^2 \sigma \ssss a_{ij}
(\ppp_id)\ppp_j(w^2) dxdtdad\tau \right\vert \\
=& \left\vert -\int_D s^2\lambda^3 \ssss \ppp_j
(\va^2\sigma a_{ij}\ppp_id)w^2 dxdtdad\tau\right\vert
\end{align*}
$$
\le C\int_D s^2\la^4\va^2 w^2 dxdtdad\tau.                     \eqno{(2.18)}
$$
\begin{align*}
&I_3 = - \int_D s\la^2\ssss \va\sigma a_{ij}w(\ppp_i\ppp_jw) dxdtdad\tau\\
=& \int_D s\la^2 \ssss \va \sigma a_{ij}(\ppp_iw)(\ppp_jw) dxdtdad\tau
+ \int_D s\la^2\ssss \ppp_i(\va\sigma a_{ij}) w (\ppp_jw)dxdtdad\tau\\
\end{align*}
$$
\ge \int_D s\la^2\va \sigma \ssss a_{ij}(\ppp_iw)(\ppp_jw) dxdtdad\tau
- C\int_D s\la^3\va\vert \nabla w\vert\vert w \vert dxdtdad\tau.   
\eqno{(2.19)}
$$
$$
I_4 = -\int_D s^3\la^4\va^3 \sigma^2 w^2 dxdtdad\tau.           \eqno{(2.20)}
$$
$$
\vert I_5 \vert \le C\left\vert \int_D s\la^2\va \times
s\la^2\va \sigma w^2 dxdtdad\tau\right\vert 
\le C\int_D s^2\la^4\va^2 w^2 dxdtdad\tau.                     \eqno{(2.21)}
$$
Hence, choosing $s>0$ and $\lambda>0$ large, by (2.16) - (2.21) we obtain
\begin{align*}
& \int_D s\la^2\va\sigma \ssss a_{ij}(\ppp_iw)(\ppp_jw) dxdtdad\tau
- \int_D s^3\la^4\va^3 \sigma^2 w^2 dxdtdad\tau\\
\le & C\int_D \vert f e^{s\va}s\la^2\va\sigma w\vert dxdtdad\tau
+ C\int_D s^2\la^4\va^2 w^2 dxdtdad\tau
+ C\int_D s\la^3\va \vert \nabla w\vert \vert w\vert dxdtdad\tau
\end{align*}
$$
\le C\int_D  f^2 \weight dxdtdad\tau + C\int_D s^2\la^4\va^2 w^2 dxdtdad\tau
+ C\int_D \la^2 \vert \nabla w\vert^2 dxdtdad\tau.      \eqno{(2.22)}
$$
At the last inequality, we argue as follows: By  
$$
s\la^3\va\vert \nabla w\vert\vert w\vert 
= (s\la^2\va \vert w\vert)(\la\vert \nabla w\vert)
\le \frac{1}{2}s^2\la^4\va^2 w^2 + \frac{1}{2}\la^2
\vert \nabla w\vert^2,
$$
we have
$$
\int_D s\la^3\va \vert \nabla w\vert\vert w\vert dxdtdad\tau
\le \frac{1}{2} \int_D (s^2\la^4\va^2 w^2 + \la^2\vert \nabla w\vert^2)
dxdtdad\tau.
$$
Furthermore
\begin{align*}
&\vert fe^{s\va}s\la^2\va \sigma w\vert \\
\le& \frac{1}{2}f^2\weight + \frac{1}{2}s^2\la^4\va^2\sigma^2w^2
\le \frac{1}{2}f^2\weight + Cs^2\la^4\va^2w^2.
\end{align*}
 
Finally we consider $2 \times (2.22) + (2.15)$.  Using
the uniform ellipticity and $\sigma_0 \equiv \inf_{x\in \Omega} 
\sigma(x) > 0$, we obtain
\begin{align*}
& \int_D s^3\la^4\va^3 \sigma_0^2 w^2 dxdtdad\tau 
+ (\sigma_0\sigma_1-C\ep)\int_D s\la^2\va\vert \nabla w\vert^2 dxdtdad\tau\\
+ &\left( \ep - \frac{C}{\la}\right) \int_D \frac{1}{s\va}
\vert \LLL w\vert^2 dxdtdad\tau \\
\le& C\int_D f^2\weight dxdtdad\tau 
\end{align*}
$$
+ C\int_D (s\la\va + \la^2)\vert \nabla w\vert^2 dxdtdad\tau
+ C\int_D (s^3\la^3\va^3 + s^2\la^4\va^2) w^2 dxdtdad\tau.       \eqno{(2.23)}
$$
Therefore, first choosing $\ep > 0$ sufficiently small such that
$\sigma_0\sigma_1 - C\ep > 0$ and then taking $\la>0$ sufficiently large
such that $\ep - \frac{C}{\la} > 0$, we can absorb the second and 
the third terms on the right-hand side of (2.23) into the left-hand 
side and we obtain
$$
\int_D s^3\la^4\va^3 w^2 dxdtdad\tau 
+ \int_D s\la^2\va\vert \nabla w\vert^2 dxdtdad\tau
+  \int_D \frac{1}{s\va} \vert \LLL w\vert^2 dxdtdad\tau 
$$
$$
\le C\int_D f^2\weight dxdtdad\tau.                              \eqno{(2.24)}
$$
Noting $w = ue^{s\va}$, we have
\begin{align*}
&\int_D \left( \frac{1}{s\va}\vert \LLL u\vert^2 
+ s\la^2\va \vert \nabla u\vert^2 + 
s^3\la^4\va^3 u^2\right) \weight dxdtdad\tau \\
\le &C\int_D f^2\weight dxdtdad\tau.           
\end{align*}
Thus the proof of Lemma 1 is completed.
\section{Proof of Theorem 1}

We need a special weight function.
The existence of such a function is proved in Fursikov and Imanuvilov
\cite{ImaF}, Imanuvilov \cite{Ima}, Imanuvilov, Puel and Yamamoto 
\cite{ImaPuY}.
\\
{\bf Lemma 2.}
Let $\omega$ be an arbitrarily fixed sub-domain of $\Omega$.
Then there exists
a fucntion $d \in C^2(\oooo\Omega)$ such that
$$
d(x) > 0 \quad x \in \Omega, \quad d\vert_{\ppp\Omega}
= 0, \quad \vert \nabla d(x)\vert > 0, \quad 
x\in \oooo{\Omega\setminus\omega}.                  
$$

{\bf Example:}
Let $\Omega = \{ x; \thinspace \vert x\vert < 1\}$ and let 
$0 \in \omega$.  Then $d(x) = 1 - \vert x\vert^2$ satisfies the conditions
in Lemma 2.

Henceforth we set 
$$
B(\p,r) := \{{\mathbf{x}} \in \R^3;\thinspace 
\vert {\mathbf{x}} - \p\vert < r\}
$$
with $\p \in \R^3$ and $r>0$, and
$$
\Vert u\Vert_{H^{1,0}(D)} = (\Vert \nabla u\Vert_{L^2(D)}^2
+ \Vert u\Vert_{L^2(D)}^2)^{\frac{1}{2}}.
$$
\\

Now we proceed to the proof of Theorem 1, which is similar to 
Theorem 5.1 in \cite{Y}.

Let $\Omega_0$ be an arbitrary subdomain of $\Omega$ such that
$\oooo \Omega_0 \subset \Omega \cup \Gamma$, $\ppp\Omega_0
\cap \ppp\Omega$ is a non-empty open subset of $\ppp\Omega$ and
$\ppp\Omega_0 \cap \ppp\Omega \subsetneqq \Gamma$, 
According to the geometry of $\Omega_0$ and
$\Gamma$, we have to choose a suitable weight function $\va$,
that is, $d(x)$.  For this, we first choose a bouned domain $\Omega_1$ 
with smooth boundary such that 
$$
\Omega \subsetneqq \Omega_1, \quad
\oooo{\Gamma} = \oooo{\ppp\Omega \cap \Omega_1}, \quad
\ppp\Omega \setminus \Gamma \subset \ppp\Omega_1.      \eqno{(3.1)}
$$
We note that $\Omega_1$ is constructed by taking a union of 
$\Omega$ and a domain $\wwww{\Omega}$ such that 
$\ppp\wwww\Omega \cap \overline{\Omega} = \Gamma$ and that 
$\Omega_1\setminus \oooo\Omega$ contains 
some non-empty open set.
Choosing $\oooo\omega \subset \Omega_1 \setminus \oooo\Omega$, 
we apply Lemma 2 to obtain $d \in C^2(\oooo\Omega_1)$ satisfying
$$
d(x) > 0, \quad x \in \Omega_1, \quad
d(x) = 0, \quad x \in \ppp\Omega_1, \quad
\vert \nabla d(x)\vert > 0, \quad x \in \oooo\Omega_1 \cap \oooo\Omega.
                                                   \eqno{(3.2)}
$$
Then, since $\overline{\Omega_0} \subset \Omega_1$, we can choose 
a sufficiently large $N > 1$ such that 
$$
\{ x\in \Omega_1; \thinspace d(x) > \frac{4}{N}\Vert d\Vert
_{C(\oooo{\Omega_1})} \} \cap \oooo\Omega \supset \Omega_0.                                                                    \eqno{(3.3)}
$$
Let $\ep>0$ be an arbitrarily small number.
Moreover we choose $\beta > 0$ such that
$$
2\beta\ep^2 > \Vert d\Vert_{C(\oooo{\Omega_1})} > \beta\ep^2.     \eqno{(3.4)}
$$
We fix $t_0 \in [\sqrt{2}\ep, T-\sqrt{2}\ep]$,
$a_0 \in [\sqrt{2}\ep, a_1-\sqrt{2}\ep]$ and $\tau_0 \in 
[\tau_1+\sqrt{2}\ep, \tau_2-\sqrt{2}\ep]$ arbitrarily.  We set 
$\p = (t_0,a_0,\tau_0)$, and 
$\va(x,t,a,\tau) = e^{\la\psi(x,t,a,\tau)}$ with fixed large parameter 
$\la>0$ and
$$
\psi(x,t,a,\tau) = d(x) - \beta ((t-t_0)^2 + (a-a_0)^2 + (\tau-\tau_0)^2)
$$
and $\mu_k = \exp\left( \la\left(\frac{k}{N}\Vert d\Vert_{C(\oooo{\Omega_1})} 
- \frac{\beta\ep^2}{N}\right)\right)$, $k=1,2,3,4$, and
$$
D = \{ (x,t,a,\tau); \thinspace x \in\oooo\Omega, \quad \va(x,t,a,\tau) > 
\mu_1 \}.
$$

Then we can verify that
$$
\Omega_0 \times B\left(\p,\frac{\ep}{\sqrt{N}}\right)
\subset D 
\subset \oooo\Omega \times B(\p, \sqrt{2}\ep).                  \eqno{(3.5)}
$$
In fact, let $(x,t,a,\tau) \in\Omega_0 \times B\left(\p, \frac{\ep}{\sqrt{N}}
\right)$.  Then, by (3.3) we have 
$x \in \oooo\Omega$ and $d(x) > \frac{4}{N}
\Vert d\Vert_{C(\oooo{\Omega_1})}$, so that
$$
d(x) - \beta((t-t_0)^2 +(a-a_0)^2 + (\tau-\tau_0)^2)
> \frac{4}{N}\Vert d\Vert_{C(\oooo{\Omega_1})}
- \frac{\beta\ep^2}{N},
$$
that is, $\va(x,t,a,\tau) > \mu_4$, which implies that $(x,t,a,\tau) \in D$
by the definition.
Next let $(x,t,a,\tau) \in D$.  Then $d(x) - \beta ((t-t_0)^2
+(a-a_0)^2 + (\tau-\tau_0)^2) > \frac{1}{N}\Vert d\Vert_{C(\oooo{\Omega_1})}
- \frac{\beta \ep^2}{N}$.  Therefore
$$
\Vert d\Vert_{C(\oooo{\Omega_1})} - \frac{1}{N}
\Vert d\Vert_{C(\oooo{\Omega_1})} + \frac{\beta\ep^2}{N}
> \beta ((t-t_0)^2+(a-a_0)^2 + (\tau-\tau_0)^2).
$$
Applying (3.4), we have $2\left(1-\frac{1}{N}\right)\beta\ep^2
+ \frac{\beta\ep^2}{N} > \beta((t-t_0)^2+(a-a_0)^2 + (\tau-\tau_0)^2)$, 
that is,
$2\beta \ep^2 > \beta((t-t_0)^2+(a-a_0)^2 + (\tau-\tau_0)^2)$.
The verification of (3.5) is completed.

Next we have 
$$
\ppp D \subset \Sigma_1 \cup \Sigma_2, \quad \mbox{where 
$\Sigma_1 \subset \Gamma \times (0,T) \times (0,a_1)\times (\tau_1,\tau_2)$},
$$
$$ 
\Sigma_2 = \{ (x,t,a,\tau); \thinspace x\in\Omega, \thinspace
\va(x,t,a,\tau) = \mu_1\}.                             \eqno{(3.6)}
$$
In fact, let $(x,t,a,\tau) \in \ppp D$.  Then $x \in \oooo\Omega$ and
$\va(x,t,a,\tau) \ge \mu_1$.  We separately consider the cases
$x \in \Omega$ and $x \in \ppp\Omega$.  First let $x \in \Omega$.
If $\va(x,t,a,\tau) > \mu_1$, then $(x,t,a,\tau)$ is an interior point of
$D$.  This is impossible.  Therefore if $x \in \Omega$, then 
$\va(x,t,a,\tau) = \mu_1$ must hold.
Next let $x \in \ppp\Omega$.  Let $x \in \ppp\Omega \setminus
\Gamma$.  Then $x \in \ppp\Omega_1$ by the third condition in 
(3.1), and $d(x) = 0$ by the second condition in (3.2).
On the other hand, $\va(x,t,a,\tau) \ge \mu_1$ yields that 
$$
d(x) - \beta((t-t_0)^2 +(a-a_0)^2 + (\tau-\tau_0)^2)) 
= -\beta((t-t_0)^2+(a-a_0)^2 + (\tau-\tau_0)^2)
\ge \frac{1}{N}\Vert d\Vert_{C(\oooo{\Omega_1})}
- \frac{\beta\ep^2}{N},
$$
that is, $0\le \beta((t-t_0)^2 +(a-a_0)^2 + (\tau-\tau_0)^2)
\le \frac{1}{N}(-\Vert d\Vert_{C(\oooo{\Omega_1})} + \beta\ep^2)$, 
which is impossible
by (3.4).  Therefore $x \in \Gamma$.  In terms of (3.5), the 
verification of (3.6) is completed.

By replacing the coefficient $c(x,t,a,\tau)$ by 
$c(x,t,a,\tau) - g'(\tau)$, it is sufficient to consider
$\LLL u - Ku = 0$ in place of $L_0u - Ku=0$.
We apply Lemma 1 in $D$ with suitably fixed $\la > 0$.
Henceforth $C>0$ denotes generic constants depending on $\la$, 
but independent of $s$.  
For it, we need a cut-off function because we have no data on 
$\ppp D \setminus (\Gamma \times (0,T)\times (0,a_1)\times (\tau_1,\tau_2))$.  
Let $\chi \in C^{\infty}(\Bbb R^{n+3})$ such that $0 \le \chi \le 1$ and
$$
\chi(x,t,a,\tau) =
\left\{
\begin{array}{rl}
1, \qquad & \va(x,t,a,\tau) > \mu_3, \\
0, \qquad & \va(x,t,a,\tau) < \mu_2.               
\end{array}
\right.                             \eqno{(3.7)}
$$

We set $v = \chi u$, and have
\begin{align*}
&Lv = (\LLL \chi) u - 2\ssss a_{ij}(\ppp_i\chi)\ppp_j u\\
-& \left( \ssss a_{ij}\ppp_i\ppp_j\chi\right)u
- \left( \sum_{i=1}^n b_i\ppp_i\chi\right) u \quad \mbox{in $D$}.
\end{align*}
Here we recall that $\LLL \chi = \ppp_t\chi + \ppp_a\chi
+ g(\tau)\ppp_\tau\chi$.
By (3.6) and (3.7), we see that
$$
v = \vert \nabla_{x,t,a,\tau} v\vert = 0 \qquad \mbox{on $\ppp D$}.
$$
Hence Lemma 1 yields
$$
\int_D s^3\vert v\vert^2 \weight dxdtda d\tau
$$
$$
\le C\int_D \left\vert
(\widetilde{L_0}\chi)u - 2\ssss a_{ij}(\ppp_i\chi)\ppp_j u
- \left( \ssss a_{ij}\ppp_i\ppp_j\chi\right)u
- \left( \sum_{i=1}^n b_i\ppp_i\chi\right) u\right\vert^2 \weight dxdtdad\tau
                                       \eqno{(3.8)}
$$
for all $s \ge s_0$.
By (3.7), the second integral on the right-hand side does not vanish
only if $\mu_2 \le \va(x,t,a,\tau) \le \mu_3$ and so
\begin{align*}
&\int_D \left\vert (\widetilde{L_0}\chi)u - 2\ssss a_{ij}(\ppp_i\chi)\ppp_j u
- \left( \ssss a_{ij}\ppp_i\ppp_j\chi\right)u
- \left( \sum_{i=1}^n b_i\ppp_i\chi\right) u\right\vert^2 \weight dxdtdad\tau\\
\le& Ce^{2s\mu_3} \Vert u\Vert^2_{H^{1,0}(D)}.
\end{align*}  
By (3.3) and the definition of $D$, we can directly verify that if 
$(x,t,a,\tau) \in \Omega_0 \times B\left(\p, \frac{\ep}{\sqrt{N}}\right)$,
then $\va(x,t,a,\tau) > \mu_4$.  
Therefore, noting (3.5) and (3.7), we see that 
\begin{align*}
&\mbox{[the left-hand side of (3.8)]}\\
\ge &\int_{B\left(\p, \frac{\ep}{\sqrt{N}}\right)}\int_{\Omega_0} 
s^3\vert v\vert^2 \weight dxdtdad\tau
\ge e^{2s\mu_4}\int_{B\left(\p, \frac{\ep}{\sqrt{N}}\right)}
\int_{\Omega_0} s^3\vert u\vert^2 dxdtdad\tau.
\end{align*}
Hence (3.8) yields
$$
e^{2s\mu_4}\int_{B\left(\p, \frac{\ep}{\sqrt{N}}\right)}
\int_{\Omega_0} s^3\vert u\vert^2 dxdtdad\tau
\le Ce^{2s\mu_3} \Vert u\Vert^2_{H^{1,0}(D)}.
$$
Therefore 
$$
\int_{B\left(\p, \frac{\ep}{\sqrt{N}}\right)}
\int_{\Omega_0} s^3\vert u\vert^2 dxdtdad\tau
\le Ce^{-2s(\mu_4-\mu_3)} \Vert u\Vert^2_{H^{1,0}(D)}
$$
for all $s \ge s_0$.  Letting $s \to \infty$, we obtain
$$
u(x,t,a,\tau) = 0, \quad 
x \in \Omega_0, \thinspace
\vert t-t_0\vert^2+\vert a-a_0\vert^2+\vert \tau-\tau_0\vert^2 
< \frac{\ep^2}{N}.                        \eqno{(3.9)}
$$
Since $(t_0,a_0,\tau_0) \in [\sqrt{2}\ep, T-\sqrt{2}\ep]
\times [\sqrt{2}\ep, a_1-\sqrt{2}\ep] \times 
[\tau_1+\sqrt{2}\ep, \tau_2-\sqrt{2}\ep]$ and $\Omega_0\subset\Omega$ 
are chosen arbitrary provided that 
$\oooo \Omega_0 \subset \Omega \cup \Gamma$, $\ppp\Omega_0 \cap 
\ppp\Omega$ is a non-empty subset of $\ppp\Omega$ and
$\ppp\Omega_0 \cap \ppp\Omega \subsetneqq \Gamma$, 
equality (3.9) yields $u=0$ in 
$\Omega \times (0,T)\times (0,a_1)\times (\tau_1,\tau_2)$.
Thus the proof of Theorem 1 is completed.

\end{document}